\documentclass[twoside,12pt]{amsart}

\usepackage[T1]{fontenc}
\usepackage{amsfonts}
\usepackage{amssymb}
\usepackage{amsthm}
\usepackage{amsmath}
\usepackage{graphicx}
\usepackage{url}
\usepackage[all]{xy}
\usepackage{lscape}
\usepackage{url}

\usepackage[shortlabels]{enumitem}
\setlist{leftmargin=*}

\theoremstyle{plain}

\newcommand{\Irr}{\mathrm{Irr}}
\newcommand{\IBr}{\mathrm{IBr}}

\newcommand{\Sp}{\mbox{\rm Sp}}

\newtheorem*{theorem*}{Theorem}

\newtheorem*{define*}{Definition}

\newtheorem*{thm*}{Theorem}

\newtheorem*{lem*}{Lemma}

\newtheorem*{prp*}{Proposition}

\newtheorem*{cor*}{Corollary}

\newtheorem*{conj*}{Conjecture}

\newtheorem*{xmpl*}{Example}

\newtheorem*{rem*}{Remark}

\begin{document}
\date{\today}
\title{The modular characters of $\Sp_{10}(2)$}

\author{Gerhard Hiss}

%\address{Lehrstuhl f\"ur Algebra und Zahlentheorie, RWTH Aachen University,
%52056 Aach\-en, Germany}
%
%\email{G.H.: gerhard.hiss@math.rwth-aachen.de}

\subjclass[2010]{20C20, 20C33, 20C40}
\keywords{Modular characters, classical group $\Sp_{10}(2)$}

\maketitle

\begin{abstract}
Using computational methods, we determine the irreducible $3$- and $5$-modular 
characters of the symplectic group $\Sp_{10}(2)$ for blocks with non cyclic 
defect groups, with one ambiguity in the principal $3$-block.
\end{abstract}

\section{Introduction and results}\label{1}

\markright{INTRODUCTION}
  
The purpose of this work is the determination of the irreducible $3$- and 
$5$-modular characters of the symplectic group $G := \Sp_{10}(2)$, where we
restrict ourselves to blocks with non-cyclic defect groups. The results were 
obtained with GAP \cite{GAP4}, in particular with the GAP package {\tt moc} 
\cite{gapmoc}. The starting point of the computations are the ordinary character 
table of~$G$, together with the ordinary and $3$- and $5$-modular character 
tables of its maximal subgroups $O_{10}^+(2)$ and $O_{10}^-(2)$; all these 
tables are available in the GAP character table library; see~\cite{CTblLib}.

Here, we present proofs of the main results for characteristic $3$, which are 
independently reproducible with any suitable computer algebra system which has 
the above basic data available, for example MAGMA~\cite{magma}. 
%%%%%%%%%%%%%%%%%%%%%%%%%%%%%%%%%%%%%%%%%%%%%%%%%%%%%%%%%%%%%%%%%%%%%%%%%%%%%%%%
%%%%%%%%%%%%%%%%%%%%%%%%%%%%%%%%%%%%%%%%%%%%%%%%%%%%%%%%%%%%%%%%%%%%%%%%%%%%%%%%
%%
%% See hiss/Uni/GapMoc/MocExp/s102/s102mod3_proof
%%
%%%%%%%%%%%%%%%%%%%%%%%%%%%%%%%%%%%%%%%%%%%%%%%%%%%%%%%%%%%%%%%%%%%%%%%%%%%%%%%%
%%%%%%%%%%%%%%%%%%%%%%%%%%%%%%%%%%%%%%%%%%%%%%%%%%%%%%%%%%%%%%%%%%%%%%%%%%%%%%%%
Except for some hints at the end of our paper, we omit the proofs for 
characteristic $5$, since they are much easier. 
%%%%%%%%%%%%%%%%%%%%%%%%%%%%%%%%%%%%%%%%%%%%%%%%%%%%%%%%%%%%%%%%%%%%%%%%%%%%%%%%
%%%%%%%%%%%%%%%%%%%%%%%%%%%%%%%%%%%%%%%%%%%%%%%%%%%%%%%%%%%%%%%%%%%%%%%%%%%%%%%%
%%
%% See hiss/Uni/GapMoc/MocExp/s102/s102mod5_proof
%%
%%%%%%%%%%%%%%%%%%%%%%%%%%%%%%%%%%%%%%%%%%%%%%%%%%%%%%%%%%%%%%%%%%%%%%%%%%%%%%%%
%%%%%%%%%%%%%%%%%%%%%%%%%%%%%%%%%%%%%%%%%%%%%%%%%%%%%%%%%%%%%%%%%%%%%%%%%%%%%%%%
For specific notation and the conceptual background of the computations, e.g., 
basic sets, we refer the reader to the monograph \cite{moc}.

In characteristic~$3$, there are five blocks with non-cyclic defect groups, with 
defects $6, 4, 4, 2, 2$, respectively. We number these blocks as in GAP by 
$B_1$, $B_2$, $B_3$, $B_5$ and $B_6$. In characteristic~$5$, there are three 
blocks with non-cyclic defect groups, all with defect~$2$, numbered $B_1$, $B_3$ 
and $B_4$. The decomposition matrices in characteristic~$3$, restricted to 
particular basic sets, are given in Tables \ref{Mod3B1PS2}, \ref{Mod3B2PS2}, 
\ref{Mod3B3PS2}, \ref{Mod3B5PS1} and \ref{Mod3B6PS1}. In characteristic~$5$,
analogous matrices are presented in Tables~\ref{Mod5B1PS1}, \ref{Mod5B3PS1} 
and \ref{Mod5B4PS1}. The irreducible Brauer characters can easily be computed 
from these matrices. The results are complete, except for one ambiguity 
in the principal $3$-block~$B_1$ of $G$.

\section{The basic sets}
\label{TheBasicSets}
From now on we work in characteristic~$3$.
We order the irreducible characters of~$G$ as in GAP; this then determines
an ordering of $\Irr( B_i )$, $i = 1, 2, 3, 5, 6$.
In order to save space, the construction of the basic sets for these blocks is 
given via a GAP program. The array {\tt positionsInobs} contains, for each $B_i$,
block, a list of the positions of the basic set characters within $\Irr( B_i )$,
$ i = 1, 2, 3, 5, 6$.

\begin{verbatim}
 s102 := CharacterTable( "S10(2)" );
 s102mod3 := CharacterTableRegular( s102, 3 );
 degs := TransposedMat( Irr( s102 ) )[ 1 ];

 ind2 := Indicator( s102, 2 );;

 nr := Filtered( [1..Length( ind2 )], i -> ind2[i] = 0 );
# [  ] 

 s102b3 := PrimeBlocks( s102, 3 );;

 defects := s102b3.defect;
# [ 6, 4, 4, 1, 2, 2, 1, 1, 0, 0, 1, 0, 0, 0 ]
 Positions( defects, 4 );
# [ 2, 3 ]
 Positions( defects, 2 );
# [ 5, 6 ]

 blnrs := [];
 for j in [1,2,3,5,6] do
 blnrs[ j ] := Filtered( [1..Length( s102b3.block )], 
                                i -> s102b3.block[ i ] = j );
 od;
# "blnrs" contains the GAP numbers of the ordinary irreducible 
# characters in the blocks 1, 2, 3, 5, 6

 ordegrees := [];
 for j in [1,2,3,5,6] do
 ordegrees[ j ] := degs{blnrs[j]}; 
 od;

 positionsInobs :=
 [ 
   [ 1, 2, 3, 5, 6, 7, 8, 9, 10, 11, 12, 13, 16, 21, 22, 25, 
    26, 27, 28, 32, 34, 37, 38, 42, 43, 45, 47, 48, 53, 55, 
    57, 68, 73, 76, 93, 109 ], 
   [ 1, 2, 3, 5, 6, 7, 9, 13, 14, 19 ], 
   [ 1, 2, 3, 5, 6, 7, 8, 12, 13, 18 ],, 
   [ 1, 2, 3, 5, 7 ], 
   [ 1, 2, 4, 5 ] 
 ];
# "positionsInobs" contains the positions of the basic set 
# characters for each block, relative to the numbering of
# the ordinary characters in the block

 restrictionsordinaries := [];
 for j in [1,2,3,5,6] do
 restrictionsordinaries[ j ] := 
 RestrictedClassFunctions( s102, Irr( s102 ){blnrs[j]}, s102mod3 ); 
 od;
# "restrictionsordinaries" contains the restrictions of the
# ordinary irreducible characters in the respective blocks to
# the 3-regular classes.

 basicsets := [];
 for j in [1,2,3,5,6] do
 basicsets[ j ] := restrictionsordinaries[j]{positionsInobs[j]};
 od;
# "basicsets" contains subsets of these, selected according to
# the positions given in "positionsInobs".

 relations := [];
 for j in [1,2,3,5,6] do
 relations[ j ] := 
     Decomposition( basicsets[j], restrictionsordinaries[j], 20 );
 od;
 
# The fact that "relations" consists of integer matrices 
# shows that the characters on Positions "positionsInobs" 
# form basic sets of ordinary characters.
\end{verbatim}

\section{The proofs}
\label{Stage1}
Except for the final subsection, we continue to work in characteristic~$3$. Let 
$B$ be one of $B_1$, $B_2$, $B_3$, $B_5$ or $B_6$, and let~$\text{BS}$ be the 
basic set of ordinary characters constructed in Section~\ref{TheBasicSets}. A 
projective character for~$B$ is given as a column, indexed by $\text{BS}$, 
recording the scalar product of the projective with the elements of $\text{BS}$. 
This information determines the projective uniquely. A \textit{relation} is a 
column of integers, indexed by $\text{BS}$, 
recording the expansion of a Brauer character of~$B$ in terms of $\text{BS}$. 
The standard inner product between a relation corresponding to the Brauer 
character $\vartheta$ and a column representing a projective character~$\Psi$, 
is the usual scalar product between $\vartheta$ and the restriction of~$\Psi$ to 
the $3$-regular classes of~$G$.

In a first phase of the proof, we construct a basic sets $\text{PS}$ of 
projective characters for~$B$. The matrix of scalar products between $\text{BS}$ 
and $\text{PS}$ is square, and it turns out that it is lower unitriangular. This 
yields an unambiguous labelling of $\IBr( B )$, and thus also of the PIMs of~$B$, 
in terms of $\text{BS}$. In a second phase, using additional projectives, we 
determine a new basic set of
projectives and show that its members (except for one case in $B_1$) are in
fact PIMs. This is achieved using relations, using the fact that the scalar 
product of a projective character and a Brauer character is a non-negative 
integer. Most of our relations arise from restrictions of ordinary characters
of $B$ to the $3$-regular classes. The projective characters originate from 
inducing PIMs from the two maximal subgroups $M^+ := O_{10}^+(2)$ and 
$M^- := O_{10}^-(2)$. The PIMs of these groups are available in GAP through the 
command {\tt DecompositionMatrix}.
In Tables~\ref{Mod3B1PS1src}--\ref{Mod3B2356PS1src}, we document the origin of
a projective by a symbol $i^{\varepsilon}$ with $i$ a positive integer and
$\varepsilon \in \{ +, - \}$. This designates the projective of~$G$ which is
induced from the $i$th projective (corresponding to the $i$th column of the
decomposition matrix of $M^{\varepsilon}$, as determined by GAP), of
$M^{\varepsilon}$ to $G$.

We will adopt the following notational convention. Projective characters of $B$ 
are denoted by $\Psi_i$, and the PIMs of $B$ by 
$\Phi_i$, $1 \leq i \leq n$, with $n = |\IBr(B)|$. We chose the notation such 
that $\text{PS} = \{ \Psi_i \mid 1 \leq i \leq n \}$. In particular, if 
$\Psi_i \in \text{PS}$ is a PIM, we have $\Psi_i = \Phi_i$. Due to the 
triangular shape of the matrix of scalar products, we have 
$\Psi_i = \Phi_i + \sum_{j = i + 1}^n a_j \Phi_j$ with non-negative integers 
$a_j$, $i + 1 \leq j \leq n$. A projective character $\Psi_i \in \text{PS}$ is 
\textit{atomic}, if the corresponding column is a unit vector. An atomic
projective is a PIM.

Occasionally, we label the irreducible ordinary or Brauer characters by their 
degrees, distinguishing those of equal degrees by subscripts, if necessary.
Consequently, the rows of the tables of scalar products and tables of relations 
given in the appendix, are labelled by the degrees of the elements of 
$\text{BS}$. As customary, zeros are replaced by dots in these tables.

\subsection{The principal block $B_1$}
Phase~$1$: 
We begin with the basic set $\text{PS}$ of projective characters given in
Table~\ref{Mod3B1PS1}. Clearly, $\Psi_{5}$, $\Psi_{23}$, $\Psi_{27}$, 
$\Psi_{30}$, $\Psi_{31}$, $\Psi_{33}$ and $\Psi_{36}$ are PIMs, since they are 
atomic. In addition, the relations show that the following projectives are PIMs:
$$
\begin{array}{c|ccccc} \hline\hline
\Psi & 12 & 19 & 22 & 29 & 35 \rule[- 4pt]{0pt}{ 16pt} \\ \hline
\text{rel} & 23, 52 & 74, 98 & 110 & 110 & 110 \rule[  0pt]{0pt}{14pt} \\ \hline\hline
\end{array}
$$

Phase~$2$: 
We use the extra projectives displayed in Table~\ref{Mod3B1Proj2} to improve on 
the basic set, and the relations in Table~\ref{Mod3B1Rels} to prove 
indecomposability of the $\Phi_i$ (except for $\Phi_{15}$). It is only for the 
indecomposability of $\Phi_{17}$, $\Phi_{11}$ and $\Phi_{1}$, where we have to 
use the relations $\otimes_{1,1}$, $\otimes_{1,2}$ and $\otimes_{1,3}$, which 
arise from tensor products of Brauer characters of~$B_1$, restricted to~$B_1$, 
as follows:
\begin{itemize}
\item   $\otimes_{1,1}$: $\varphi_3 \otimes \varphi_6 = 186 \otimes 10881$;
\item   $\otimes_{1,2}$: $\varphi_2 \otimes \varphi_6 = 155 \otimes 10881$;
\item   $\otimes_{1,3}$: $\varphi_3 \otimes \varphi_{17} = 186 \otimes 133672$.
\end{itemize}
Here, $\varphi_i$ denotes the irreducible Brauer character corresponding to the
PIM $\Phi_i$.
We emphasize that the existence of these Brauer characters is guaranteed at the 
stages where they are employed. The proof is given in note form. If a new 
projective is atomic, we suppress the remark that it is a PIM.

%curPS11expBSPhase2 := [];
%curPS11expBSPhase2[19] := curPS11expBS[19];
%curPS11expBSPhase2[20] := curPS11expBS[20] - curPS11expBS[33]; # wegen 55
%curPS11expBSPhase2[21] := curPS11expBS[21] - curPS11expBS[36]; # wegen 67
%curPS11expBSPhase2[22] := curPS11expBS[22];
%curPS11expBSPhase2[23] := curPS11expBS[23];
%curPS11expBSPhase2[24] := curPS11expBS[24] - 2*curPS11expBS[32] + curPS11expBS[36];
%                                             # wegen 50 und Relation 70
%curPS11expBSPhase2[25] := curPS11expBS[25] - curPS11expBS[36]; # wegen 51
%curPS11expBSPhase2[26] := curPS11expBS[26] - 2*curPS11expBS[35]; # wegen 42
%curPS11expBSPhase2[27] := curPS11expBS[27];
%curPS11expBSPhase2[28] := curPS11expBS[28] - 2 * curPS11expBS[35]; # wegen 39 und 61
%curPS11expBSPhase2[29] := curPS11expBS[29];
%curPS11expBSPhase2[30] := curPS11expBS[30];
%curPS11expBSPhase2[31] := curPS11expBS[31];
%curPS11expBSPhase2[32] := curPS11expBS[32] - 2 * curPS11expBS[36]; # wegen 64
%curPS11expBSPhase2[33] := curPS11expBS[33];
%curPS11expBSPhase2[34] := curPS11expBS[34] - curPS11expBS[35]; # wegen 37
%curPS11expBSPhase2[35] := curPS11expBS[35];
%curPS11expBSPhase2[36] := curPS11expBS[36];
\begin{enumerate}
\item $\Phi_{34} := \Psi_{34} - \Phi_{35}$ is projective by $\Psi_{37}$;
\item $\Psi_{39}$ contains contains $\Phi_{24}$ three times. In particular, the 
      multiplicity of $2213400$ as a constituent of $\Phi_{24}$ is at most~$1$.
      It follows that $\Phi_{32}$ is contained at least twice in $\Psi_{24}$.
      Hence $\Phi_{32} = \Psi_{32} - 2\cdot \Phi_{36}$;
\item $\Psi_{43}$ contains $\Phi_{28}$ twice, and so $\Phi_{28} - 
      2 \cdot \Phi_{35}$ is projective; it is in fact a PIM by relation~$98$;
\item $\Psi_{38}$ contains $\Phi_{26}$ twice, and so $\Phi_{26} = \Psi_{26} -
      2 \cdot \Phi_{35}$;
\item $\Phi_{25} := \Psi_{25} - \Phi_{36}$ is projective by $\Psi_{40}$, 
      and a PIM by relation~$98$;
\item $\Psi_{39}$ contains $\Phi_{24}$ three times. Also, $\Phi_{35}$ must be 
      contained in $\Psi_{39}$. It follows that $18260550$ is not a constituent 
      of $\Phi_{24}$. By relation $70$, the PIM $\Phi_{24}$ must contain the 
      ordinary character $2213400$. This yields 
      $\Phi_{24} = \Psi_{24} - 2\cdot \Phi_{32} - 3 \cdot \Phi_{36}$;
\item $\Phi_{21} := \Psi_{21} - \Phi_{36}$ is projective by $\Psi_{46}$, 
      and a PIM by relation~$77$; 
\item $\Phi_{20} := \Psi_{20} - \Phi_{33}$ is projective by $\Psi_{41}$, 
      and a PIM by relation~$98$;
\item $\Phi_{18} := \Psi_{18} - 2 \cdot \Phi_{23}  - 2 \cdot \Phi_{25}$
      is projective by $\Psi_{45}$, and a PIM by relations~$51$,~$60$ and~$98$;
\item $\Psi_{17} := \Psi_{17}$ is a PIM by relations~$33$, $49$, $59$, $104$ 
      and~$\otimes_{1,3}$;
\item $\Phi_{16} := \Psi_{16} - 2 \cdot \Phi_{34} - \Phi_{35}$ is projective
      by $\Psi_{42}$ and $\Psi_{43}$, and a PIM by relation~$80$;
\item $\Psi_{15}^{[2]} := \Psi_{15} - 2 \cdot \Phi_{31} - \Phi_{32}$ is 
      projective by $\Psi_{44}$; by relations~$70$ and~$112$, the characters
      $1623160_2$ and $2213400$ are constituents of $\Phi_{15}$;
\item $\Phi_{14} := \Psi_{14} - \Phi_{29} - 2 \cdot \Phi_{30}$ is projective;
      indeed, $\Phi_{29}$ is contained in $\Psi_{14}$ by~$\Psi_{46}$; also,
      $\Psi_{52}$ contains $\Phi_{14}$ at least twice, and so the multiplicity
      of $1623160_1$ in $\Phi_{14}$ is at most~$1$; relations~$80$ and~$98$
      show that $\Psi_{14} - \Phi_{29} - 2 \cdot \Phi_{30}$ is a PIM;
\item $\Phi_{13} := \Psi_{13} - 2 \cdot \Phi_{19} - \Phi_{28}$ is 
      projective by $\Psi_{47}$; it is a PIM by relations $29$, $30$, $40$, $66$, 
      $84$ and~$88$;
\item $\Psi_{11} := \Psi_{11} - \Phi_{34} - \Phi_{35}$ is projective by 
      $\Psi_{48}$; it is a PIM by relations~$54$, $80$ and~$\otimes_{1,2}$;
\item $\Phi_{10} := \Psi_{10} - \Phi_{26} - 2 \dot \Phi_{32} - \Phi_{36}$ 
      is  projective by $\Psi_{49}$; indeed, $\Phi_{15}$ is contained twice in 
      $\Psi_{49}$, and $2213400$ is a constituent of $\Phi_{15}$, so that the
      multiplicity of $2213400$ in $\Phi_{10}$ is at most $1$; also, 
      $\Psi_{10} - \Phi_{26} - 2 \dot \Phi_{32} - \Phi_{36}$ is a PIM by 
      relations~$23$ and~$110$;
\item $\Phi_{9} := \Psi_{9} - \Phi_{35}$ is projective by $\Psi_{50}$,
      and a PIM by relation~$40$,~$56$,~$71$;
\item $\Phi_{8} := \Psi_{8} - 2 \cdot \Phi_{26}$ is projective by 
      $\Psi_{50}$, and the fact that $\Phi_9$ contains $973896$;
\item $\Phi_{7} := \Psi_{7} - \Phi_{30}$ is projective by $\Psi_{51}$ 
      and a PIM by relations~$23$ and~$105$;
\item $\Phi_{6} := \Psi_{6} - 2 \cdot \Phi_{20} - 2 \cdot \Phi_{23}$ is 
      projective by $\Psi_{53}$, and a PIM by relations~$23$ and~$65$;
\item $\Phi_4 := 527 + 11594_2 + 17050$ is projective by $\Psi_4, \Psi_{54},
      \Psi_{55}$ and $\Psi_{56}$, and a PIM by relations~$14$ and~$29$;
\item $\Phi_3 := \Psi_3 - 2 \cdot \Phi_{23} - \Phi_{25}$
      is projective by $\Psi_{57}$ and $\Psi_{58}$, and a PIM by relations~$18$,
      $19$, $110$, $40$ and $65$;
\item $\Psi_2 := (\Psi_{59} - \Phi_4 - \Phi_{12})/2 - \Phi_{22}$
      is projective by $\Psi_2$, $\Psi_{59}$ and $\Psi_{60}$; $\Psi_2$ is a PIM 
      by relations~$17$, $77$, $79$ and $82$ and the fact that neither 
      $\Phi_{13}$ nor $\Phi_{20}$ are components of $\Psi_2$;
\item $\Psi_1 := \Psi_1 - \Phi_5 - 2 \cdot \Phi_7 - \Phi_8 - \Phi_{27} - \Phi_{30}$
      is projective by $\Psi_{61}$; it is a PIM by relations~$4$ 
      and $\otimes_{1,1}$.
\end{enumerate}
This proves the correctness of Table~\ref{Mod3B1PS2}.

\subsection{The block $B_2$}
A first basic set of projective characters for $B_2$, as well as some further
projectives are given in Table~\ref{Mod3B2PS1}. The origin of these projectives 
is recorded in Table~\ref{Mod3B2356PS1src}. Notice that in particular $\Psi_6$ 
and $\Psi_{10}$ are PIMs, since they are atomic.
Projective $\Psi_{11}$ shows that $\Psi_1 - \Phi_6$ is projective.
Projective $\Psi_{12}$ shows that $\Psi_3 - \Phi_{10}$ is projective.
Projective $\Psi_{13}$ shows that $\Psi_7 - \Phi_{10}$ is projective.
Projective $\Psi_{14}$ shows that $\Psi_8 - \Phi_{10}$ is projective.

This yields the second basic set of projectives displayed in
Table~\ref{Mod3B2PS2}. This table also contains the relations arising from the 
characters $\chi_{4}$, $\chi_{10}$, $\chi_{11}$, $\chi_{16}$, $\chi_{20}$ and 
$\chi_{23}$ (numbering relative to $\Irr( B_2 )$); the relation headed 
$\otimes_{2,1}$ arises from the product $\varphi_2 \otimes \varphi_{20} = 
155 \otimes 237181$, restricted to $B_2$, where $\varphi_2$ and $\varphi_{20}$ 
are the irreducible Brauer characters of $B_1$ corresponding to the PIMs 
$\Phi_2$ and $\Phi_{20}$, respectively. Relations~$10$ and~$16$ show that 
$\Phi_2$ is a PIM. Relations~$11$ and~$20$ show that $\Phi_3$ is a PIM.
Relation~$23$ shows that $\Psi_4$ and $\Psi_8$ are PIMs. Finally, relations~$4$
and~$\otimes_{2,1}$ show that $\Psi_1$ is a PIM.
This proves the correctness of Table~\ref{Mod3B2PS2}.

\subsection{The block $B_3$}
A first basic set of projective characters for $B_3$, as well as one further
projective are given in Table~\ref{Mod3B3PS1}. The origin of these projectives 
is given in Table~\ref{Mod3B2356PS1src}. Table~\ref{Mod3B3PS1} also contains the 
relations arising from the characters $\chi_{15}$, $\chi_{22}$ and $\chi_{20}$ 
(numbering relative to $\Irr( B_3 )$); the relation headed $\otimes_{3,1}$ 
arises from the product $\varphi_2 \otimes \varphi_{8} = 155 \otimes 16337$, 
restricted to $B_3$, where $\varphi_2$ and $\varphi_{8}$ are the irreducible 
Brauer characters of~$B_1$ corresponding to the PIMs $\Phi_2$ and $\Phi_{8}$, 
respectively.  Relations~$15$, $20$ and $22$ show that $\Phi_2$, $\Phi_3$ and 
$\Phi_7$ are indecomposable. The projective $\Phi_{11}$ shows that 
$\Phi_5 - \Phi_9 - 2\Phi_{10}$ is projective. Relation~$20$ shows that the 
latter is a PIM. This yields the second basic set of projectives given in 
Table~\ref{Mod3B3PS2}. Finally, the tensor relations shows that $\Phi_1$ is a 
PIM, since $\Phi_3$ cannot be contained in $\Psi_1$.
This proves the correctness of Table~\ref{Mod3B3PS2}.

\subsection{The block $B_5$}
We begin with the basic set of projective characters for $B_5$, given in
Table~\ref{Mod3B5PS1}. The origin of these projectives is given in
Table~\ref{Mod3B2356PS1src}. The relation arising from character $\chi_{8}$
(numbering relative to $\Irr( B_3 )$) is given in the last column of 
Table~\ref{Mod3B5PS1}. It follows from this relation that $\Phi_2$ and $\Phi_3$ 
are PIMs. In turn, $\Phi_1$ is a PIM, since $\Phi_3$ cannot be contained in 
$\Phi_1$. This proves the correctness of Table~\ref{Mod3B5PS1}.

\subsection{The block $B_6$}
We begin with the basic set of projective characters for $B_6$, given in
Table~\ref{Mod3B6PS1}. The origin of these projectives is given in
Table~\ref{Mod3B2356PS1src}. These projectives are evidently indecomposable,
proving our result for $B_6$.

\subsection{On the proofs in characteristic $5$}
In characteristic~$5$, all but one of the PIMs of the blocks $B_1$, $B_3$ and 
$B_4$ are induced from PIMs of $O_{10}^+(2)$; the exception is the PIM 
$\Phi_7$ of $B_1$, which is induced from a PIM of $O_{10}^-(2)$. The
indecomposability of the given projectives is easily deduced from relations
and the triangularity of the decomposition matrices.

\bigskip
\bigskip

{\small
\noindent\textsc{Lehrstuhl f\"ur Algebra und Zahlentheorie, RWTH Aachen University,
52056 Aach\-en, Germany}

\noindent\textit{Email address: }\texttt{gerhard.hiss@math.rwth-aachen.de}
}

\newpage

%\cleardoublepage

\clearpage

\appendix

\section{Tables}

\begin{table}[h]
\caption{\label{Mod3B1PS1} The first basic set of projectives for $B_1$} 
$$
{\footnotesize
\begin{array}{r|rrrrrrrrrrrrrrrrrr} \hline\hline
\mbox{PIM}&   &   &   &   & * &   &   &   &   &   &   & * &   &   &   &   &   &  \\ \hline
          &   &   &   &   &   &   &   &   &   & 1 & 1 & 1 & 1 & 1 & 1 & 1 & 1 & 1 \\
\Psi      & 1 & 2 & 3 & 4 & 5 & 6 & 7 & 8 & 9 & 0 & 1 & 2 & 3 & 4 & 5 & 6 & 7 & 8 \rule[- 4pt]{0pt}{ 14pt} \\ \hline     
        1 & 1 & . & . & . & . & . & . & . & . & . & . & . & . & . & . & . & . & . \\
      155 & . & 1 & . & . & . & . & . & . & . & . & . & . & . & . & . & . & . & . \\
      187 & 1 & . & 1 & . & . & . & . & . & . & . & . & . & . & . & . & . & . & . \\
      527 & 1 & . & . & 1 & . & . & . & . & . & . & . & . & . & . & . & . & . & . \\
    11594 & 1 & . & . & . & 1 & . & . & . & . & . & . & . & . & . & . & . & . & . \\
    11594 & 1 & . & 1 & 1 & . & 1 & . & . & . & . & . & . & . & . & . & . & . & . \\
    11935 & 2 & . & . & . & . & . & 1 & . & . & . & . & . & . & . & . & . & . & . \\
    17050 & 2 & 1 & 1 & 1 & . & . & . & 1 & . & . & . & . & . & . & . & . & . & . \\
    17391 & . & . & 1 & . & . & . & . & . & 1 & . & . & . & . & . & . & . & . & . \\
    18445 & . & . & . & . & . & . & . & . & . & 1 & . & . & . & . & . & . & . & . \\
    28985 & . & 1 & . & 1 & . & . & . & . & . & . & 1 & . & . & . & . & . & . & . \\
    29512 & . & 3 & . & . & . & . & . & . & . & . & . & 1 & . & . & . & . & . & . \\
    46376 & . & 1 & . & . & . & . & . & . & . & . & . & . & 1 & . & . & . & . & . \\
    81158 & 2 & . & . & . & . & . & 1 & . & . & . & . & . & . & 1 & . & . & . & . \\
   173910 & . & 2 & . & . & . & 1 & . & . & . & 1 & . & 1 & . & . & 1 & . & . & . \\
   202895 & . & . & . & . & . & . & . & . & . & . & . & . & 1 & . & . & 1 & . & . \\
   202895 & . & . & . & . & . & . & . & . & . & . & . & . & . & 1 & . & . & 1 & . \\
   202895 & . & . & . & . & . & . & . & . & . & . & . & . & 1 & . & . & . & . & 1 \\
   202895 & . & . & . & . & . & . & . & . & . & . & . & . & 3 & . & . & . & . & . \\
   371008 & . & 1 & . & . & . & 2 & . & . & . & . & . & . & . & . & . & . & 1 & . \\
   405790 & . & . & . & 1 & . & . & . & . & . & . & . & . & . & . & . & . & . & . \\
   405790 & . & 1 & . & . & . & . & . & . & . & . & . & . & 1 & . & . & . & . & . \\
   405790 & . & 1 & 3 & . & . & 2 & . & . & 1 & . & . & . & 1 & . & . & . & . & 2 \\
   811580 & . & . & . & . & . & . & . & . & . & . & . & . & . & . & . & . & 1 & . \\
   869550 & . & . & 1 & 2 & . & . & . & . & 1 & . & 1 & . & . & . & . & . & . & 3 \\
   973896 & . & . & 1 & 2 & . & 1 & . & 2 & 1 & 1 & . & . & . & . & . & . & . & . \\
  1014475 & 1 & 3 & . & 3 & . & . & . & . & . & . & 1 & 1 & . & . & . & . & . & . \\
  1014475 & . & . & . & . & . & . & . & . & . & . & . & . & 1 & . & . & . & . & 1 \\
  1298528 & 2 & . & . & . & . & . & 1 & . & . & . & . & . & . & 1 & . & . & . & . \\
  1623160 & 1 & . & . & . & . & . & 1 & . & . & . & . & . & . & 3 & . & . & 1 & . \\
  1623160 & . & . & . & . & . & . & . & . & . & . & . & . & 3 & . & 3 & . & . & . \\
  2213400 & . & . & . & . & . & . & . & . & . & 3 & . & . & . & . & 2 & . & . & . \\
  3246320 & . & 1 & . & . & . & . & . & . & . & . & . & . & 1 & . & 1 & . & . & . \\
  3478200 & . & 1 & . & 2 & . & . & . & . & . & . & 2 & . & 1 & 1 & . & 3 & . & . \\
  8115800 & . & . & 1 & 1 & . & 2 & . & . & 1 & . & 1 & . & 5 & . & . & 1 & 1 & 3 \\
 18260550 & . & . & . & 1 & . & . & . & . & 1 & 1 & 1 & . & . & 1 & . & 1 & 2 & . 
            \rule[ 0pt]{0pt}{ 4pt} \\ \hline\hline
\end{array}
}
$$
\end{table}

\clearpage

\addtocounter{table}{-1}
\begin{table}[h]
\caption{The first basic set of projectives for $B_1$ (cont)}
$$
{\footnotesize
\begin{array}{r|rrrrrrrrrrrrrrrrrr|r} \hline\hline
\mbox{PIM}& * &   &   & * & * &   &   &   & * &   & * & * & * &   & * &   & * & * &  \\ \hline
          & 1 & 2 & 2 & 2 & 2 & 2 & 2 & 2 & 2 & 2 & 2 & 3 & 3 & 3 & 3 & 3 & 3 & 3 & 1 \\
\Psi      & 9 & 0 & 1 & 2 & 3 & 4 & 5 & 6 & 7 & 8 & 9 & 0 & 1 & 2 & 3 & 4 & 5 & 6 & 1' \rule[- 4pt]{0pt}{ 14pt} \\ \hline     
        1 & . & . & . & . & . & . & . & . & . & . & . & . & . & . & . & . & . & . & . \\
      155 & . & . & . & . & . & . & . & . & . & . & . & . & . & . & . & . & . & . & . \\
      187 & . & . & . & . & . & . & . & . & . & . & . & . & . & . & . & . & . & . & . \\
      527 & . & . & . & . & . & . & . & . & . & . & . & . & . & . & . & . & . & . & . \\
    11594 & . & . & . & . & . & . & . & . & . & . & . & . & . & . & . & . & . & . & . \\
    11594 & . & . & . & . & . & . & . & . & . & . & . & . & . & . & . & . & . & . & . \\
    11935 & . & . & . & . & . & . & . & . & . & . & . & . & . & . & . & . & . & . & . \\
    17050 & . & . & . & . & . & . & . & . & . & . & . & . & . & . & . & . & . & . & . \\
    17391 & . & . & . & . & . & . & . & . & . & . & . & . & . & . & . & . & . & . & . \\
    18445 & . & . & . & . & . & . & . & . & . & . & . & . & . & . & . & . & . & . & . \\
    28985 & . & . & . & . & . & . & . & . & . & . & . & . & . & . & . & . & . & . & 2 \\
    29512 & . & . & . & . & . & . & . & . & . & . & . & . & . & . & . & . & . & . & . \\
    46376 & . & . & . & . & . & . & . & . & . & . & . & . & . & . & . & . & . & . & . \\
    81158 & . & . & . & . & . & . & . & . & . & . & . & . & . & . & . & . & . & . & . \\
   173910 & . & . & . & . & . & . & . & . & . & . & . & . & . & . & . & . & . & . & . \\
   202895 & . & . & . & . & . & . & . & . & . & . & . & . & . & . & . & . & . & . & . \\
   202895 & . & . & . & . & . & . & . & . & . & . & . & . & . & . & . & . & . & . & . \\
   202895 & . & . & . & . & . & . & . & . & . & . & . & . & . & . & . & . & . & . & . \\
   202895 & 1 & . & . & . & . & . & . & . & . & . & . & . & . & . & . & . & . & . & . \\
   371008 & . & 1 & . & . & . & . & . & . & . & . & . & . & . & . & . & . & . & . & . \\
   405790 & . & . & 1 & . & . & . & . & . & . & . & . & . & . & . & . & . & . & . & . \\
   405790 & . & . & . & 1 & . & . & . & . & . & . & . & . & . & . & . & . & . & . & . \\
   405790 & . & . & . & . & 1 & . & . & . & . & . & . & . & . & . & . & . & . & . & . \\
   811580 & . & . & . & . & . & 1 & . & . & . & . & . & . & . & . & . & . & . & . & . \\
   869550 & . & . & . & . & . & . & 1 & . & . & . & . & . & . & . & . & . & . & . & 2 \\
   973896 & . & . & . & . & . & . & . & 1 & . & . & . & . & . & . & . & . & . & . & . \\
  1014475 & . & . & . & . & . & . & . & . & 1 & . & . & . & . & . & . & . & . & . & 2 \\
  1014475 & . & . & . & . & . & . & . & . & . & 1 & . & . & . & . & . & . & . & . & . \\
  1298528 & . & . & . & . & . & . & . & . & . & . & 1 & . & . & . & . & . & . & . & . \\
  1623160 & . & . & . & . & . & . & . & . & . & . & . & 1 & . & . & . & . & . & . & . \\
  1623160 & 1 & . & . & . & . & . & . & . & . & . & . & . & 1 & . & . & . & . & . & . \\
  2213400 & . & . & . & . & . & 3 & . & . & . & . & . & . & . & 1 & . & . & . & . & . \\
  3246320 & . & 1 & . & 1 & . & . & . & . & . & . & . & . & . & . & 1 & . & . & . & . \\
  3478200 & . & . & 1 & . & . & . & . & . & . & . & . & . & . & . & . & 1 & . & . & 3 \\
  8115800 & 2 & 1 & . & . & . & . & 1 & 2 & . & 3 & . & . & . & . & . & 1 & 1 & . & 1 \\
 18260550 & . & . & 1 & . & . & 3 & 1 & 2 & . & 2 & 1 & . & . & 2 & . & 1 & 1 & 1 & 1  
            \rule[ 0pt]{0pt}{ 4pt} \\ \hline\hline
\end{array}
}
$$
\end{table}

\clearpage

\begin{table}[h]
\caption{\label{Mod3B1Proj2} Some further projectives for $B_1$}
$$
{\small
\begin{array}{r|rrrrrrrrrrrr} \hline\hline
\Psi     & 3 & 3 & 3 & 4 & 4 & 4 & 4 & 4 & 4 & 4 & 4 & 4\\
         & 7 & 8 & 9 & 0 & 1 & 2 & 3 & 4 & 5 & 6 & 7 & 8 \rule[- 4pt]{0pt}{ 14pt} \\ \hline     
       1 & . & . & . & . & . & . & . & . & . & . & . & . \\
     155 & . & . & . & . & . & . & . & . & . & . & . & . \\
     187 & . & . & . & . & . & . & . & . & . & . & . & . \\
     527 & . & . & . & . & . & . & . & . & . & . & . & . \\
   11594 & . & . & . & . & . & . & . & . & . & . & . & 2 \\
   11594 & . & . & . & . & . & . & . & . & . & . & . & . \\
   11935 & . & . & . & . & . & . & . & . & . & . & . & . \\
   17050 & . & . & . & . & . & . & . & . & . & . & . & . \\
   17391 & . & . & . & . & . & . & . & . & . & . & . & . \\
   18445 & . & . & . & . & . & . & . & . & . & . & . & . \\
   28985 & . & . & . & . & . & . & . & . & . & . & . & 1 \\
   29512 & . & . & . & . & . & . & . & . & . & . & . & . \\
   46376 & . & . & . & . & . & . & . & . & . & . & 2 & . \\
   81158 & . & . & . & . & . & . & . & . & . & 1 & . & . \\
  173910 & . & . & . & . & . & . & . & 1 & 1 & . & . & . \\
  202895 & . & . & . & . & . & 1 & 1 & . & . & . & 2 & . \\
  202895 & . & . & . & . & . & . & . & . & . & 1 & . & . \\
  202895 & . & . & . & . & . & . & . & . & 1 & . & 2 & . \\
  202895 & . & . & . & . & . & . & . & . & . & . & 2 & . \\
  371008 & . & . & . & . & 1 & . & . & . & 2 & . & . & . \\
  405790 & . & . & . & . & . & . & . & . & . & 1 & . & . \\
  405790 & . & . & . & . & . & . & . & 1 & . & . & 4 & . \\
  405790 & . & . & . & 1 & 1 & . & 2 & . & . & . & . & . \\
  811580 & . & . & 3 & . & . & . & . & . & . & . & . & . \\
  869550 & . & . & . & 1 & 1 & . & . & . & 1 & . & . & 2 \\
  973896 & . & 2 & . & . & . & . & . & . & . & 1 & . & . \\
 1014475 & . & . & . & . & . & 2 & . & . & . & . & . & 1 \\
 1014475 & . & . & . & . & . & . & 2 & . & 1 & . & . & . \\
 1298528 & . & . & . & . & . & 2 & . & . & . & . & . & . \\
 1623160 & 1 & . & . & . & . & . & . & . & . & 3 & . & . \\
 1623160 & . & . & . & . & . & . & . & 1 & 1 & 2 & 2 & . \\
 2213400 & . & . & 4 & . & . & . & . & 1 & 1 & . & . & . \\
 3246320 & . & . & . & . & . & . & . & 4 & 3 & 1 & 4 & . \\
 3478200 & 1 & . & 2 & . & . & 2 & 1 & . & . & 2 & 6 & 1 \\
 8115800 & . & 1 & 1 & 1 & 2 & . & 3 & . & 3 & . & 1 & 1 \\
18260550 & 1 & 1 & 3 & . & . & 2 & 1 & . & . & . & 1 & . 
            \rule[ 0pt]{0pt}{ 4pt} \\ \hline\hline
\end{array}
}
$$
\end{table}

% ll5{[1..14]};
%[ 37, 42, 50, 51, 55, 60, 61, 64, 66, 67, 69, 70, 72, 77 ] 
% omitted 69, 70

\addtocounter{table}{-1}
\begin{table}[h]
\caption{Some further projectives for $B_1$ (cont)}
$$
{\small
\begin{array}{r|rrrrrrrrrrrrr} \hline\hline
\Psi     & 4 & 5 & 5 & 5 & 5 & 5 & 5 & 5 & 5 & 5 & 5 & 6 & 6 \\
         & 9 & 0 & 1 & 2 & 3 & 4 & 5 & 6 & 7 & 8 & 9 & 0 & 1 \rule[- 4pt]{0pt}{ 14pt} \\ \hline     
       1 & . & . & . & . & . & . & . & . & . & . & . & . & 1 \\
     155 & . & . & . & . & . & . & . & . & . & . & 2 & 1 & . \\
     187 & . & . & . & . & . & . & . & . & 1 & 1 & . & . & 1 \\
     527 & . & . & . & . & . & 1 & 2 & 2 & . & 1 & 1 & . & 1 \\
   11594 & . & 1 & 1 & 1 & . & 2 & . & . & . & . & . & . & . \\
   11594 & . & . & . & . & 1 & 1 & 2 & 2 & 1 & 2 & 1 & . & 1 \\
   11935 & . & . & 1 & 1 & . & 2 & . & . & . & . & . & . & . \\
   17050 & . & 1 & . & 1 & . & 1 & 2 & 2 & 1 & 2 & 1 & . & 1 \\
   17391 & . & 1 & 1 & . & . & . & . & 1 & 1 & 1 & . & 1 & 1 \\
   18445 & 1 & . & . & . & . & . & . & . & . & . & . & . & . \\
   28985 & . & . & . & . & . & . & 2 & . & . & . & 2 & 1 & . \\
   29512 & 1 & . & . & . & . & . & . & . & 1 & . & 3 & 1 & . \\
   46376 & 1 & . & . & . & . & . & . & . & . & . & 2 & 1 & . \\
   81158 & . & . & 1 & 3 & . & 2 & . & . & . & . & . & . & . \\
  173910 & 4 & . & . & . & 1 & . & . & . & 1 & . & 1 & . & . \\
  202895 & 1 & . & 1 & 2 & . & . & . & . & . & . & . & . & . \\
  202895 & . & . & . & 2 & . & . & . & . & . & . & . & . & . \\
  202895 & 1 & . & . & . & 1 & . & . & . & . & . & . & . & . \\
  202895 & 1 & . & 2 & . & . & . & . & . & . & . & . & 1 & . \\
  371008 & 1 & . & . & . & . & . & . & . & 2 & . & 2 & 1 & . \\
  405790 & . & . & . & . & . & . & . & . & . & . & . & . & . \\
  405790 & 3 & . & 1 & . & . & . & . & . & . & . & 2 & . & . \\
  405790 & 1 & 1 & 1 & . & . & . & . & 1 & 2 & 1 & 2 & 4 & 3 \\
  811580 & . & . & . & . & . & . & . & . & . & . & . & . & . \\
  869550 & . & 1 & 1 & . & 3 & . & 2 & 1 & . & 1 & . & 1 & 1 \\
  973896 & . & 1 & 1 & . & 3 & . & . & 3 & 1 & 3 & . & 1 & 1 \\
 1014475 & 1 & . & . & 2 & 1 & . & 3 & . & 1 & . & 3 & 1 & . \\
 1014475 & . & . & . & . & 1 & . & . & . & . & . & . & . & . \\
 1298528 & . & . & 3 & 3 & . & 2 & 2 & . & . & . & . & . & . \\
 1623160 & . & . & . & 3 & . & . & . & . & . & . & . & . & . \\
 1623160 & 3 & . & 2 & . & . & . & . & . & . & . & . & 1 & . \\
 2213400 & 3 & . & . & . & 1 & . & . & . & . & . & . & . & . \\
 3246320 & 6 & . & 1 & . & . & . & . & . & 1 & . & 2 & . & . \\
 3478200 & 1 & . & 3 & 5 & . & . & 2 & . & . & . & 2 & 1 & . \\
 8115800 & 1 & . & 5 & . & 3 & . & . & 1 & 2 & 1 & . & 2 & 1 \\
18260550 & . & . & 3 & 2 & 2 & . & 2 & 1 & . & 1 & . & . & 1 
            \rule[ 0pt]{0pt}{ 4pt} \\ \hline\hline
\end{array}
}
$$
\end{table}

%ll5{[15..28]};
%[ 81, 84, 86, 87, 90, 93, 94, 95, 96, 97, 99, 35, 36, 100 ]
% omitted 100

\clearpage

\begin{table}[h]
\caption{\label{Mod3B1PS2} The decomposition matrix modulo $3$ for $B_1$}
{\footnotesize
$$
\begin{array}{r|rrrrrrrrrrrrrrrrrr} \hline\hline
\mbox{PIM}& * & * & * & * & * & * & * & * & * & * & * & * & * & * & * & * & * & * \\ \hline
          &   &   &   &   &   &   &   &   &   & 1 & 1 & 1 & 1 & 1 & 1 & 1 & 1 & 1 \\
\Phi      & 1 & 2 & 3 & 4 & 5 & 6 & 7 & 8 & 9 & 0 & 1 & 2 & 3 & 4 & 5 & 6 & 7 & 8 \rule[- 4pt]{0pt}{ 14pt} \\ \hline     
       1  & 1 & . & . & . & . & . & . & . & . & . & . & . & . & . & . & . & . & . \\ 
     155  & . & 1 & . & . & . & . & . & . & . & . & . & . & . & . & . & . & . & . \\ 
     187  & 1 & . & 1 & . & . & . & . & . & . & . & . & . & . & . & . & . & . & . \\ 
     527  & 1 & . & . & 1 & . & . & . & . & . & . & . & . & . & . & . & . & . & . \\ 
   11594  & . & . & . & . & 1 & . & . & . & . & . & . & . & . & . & . & . & . & . \\ 
   11594  & 1 & . & 1 & 1 & . & 1 & . & . & . & . & . & . & . & . & . & . & . & . \\ 
   11935  & . & . & . & . & . & . & 1 & . & . & . & . & . & . & . & . & . & . & . \\ 
   17050  & 1 & . & 1 & 1 & . & . & . & 1 & . & . & . & . & . & . & . & . & . & . \\ 
   17391  & . & . & 1 & . & . & . & . & . & 1 & . & . & . & . & . & . & . & . & . \\ 
   18445  & . & . & . & . & . & . & . & . & . & 1 & . & . & . & . & . & . & . & . \\ 
   28985  & . & 1 & . & . & . & . & . & . & . & . & 1 & . & . & . & . & . & . & . \\ 
   29512  & . & 1 & . & . & . & . & . & . & . & . & . & 1 & . & . & . & . & . & . \\ 
   46376  & . & 1 & . & . & . & . & . & . & . & . & . & . & 1 & . & . & . & . & . \\ 
   81158  & . & . & . & . & . & . & 1 & . & . & . & . & . & . & 1 & . & . & . & . \\ 
  173910  & . & . & . & . & . & 1 & . & . & . & 1 & . & 1 & . & . & 1 & . & . & . \\ 
  202895  & . & . & . & . & . & . & . & . & . & . & . & . & 1 & . & . & 1 & . & . \\ 
  202895  & . & . & . & . & . & . & . & . & . & . & . & . & . & 1 & . & . & 1 & . \\ 
  202895  & . & . & . & . & . & . & . & . & . & . & . & . & 1 & . & . & . & . & 1 \\ 
  202895  & . & . & . & . & . & . & . & . & . & . & . & . & 1 & . & . & . & . & . \\ 
  371008  & . & 1 & . & . & . & . & . & . & . & . & . & . & . & . & . & . & 1 & . \\ 
  405790  & . & . & . & . & . & . & . & . & . & . & . & . & . & . & . & . & . & . \\ 
  405790  & . & . & . & . & . & . & . & . & . & . & . & . & 1 & . & . & . & . & . \\ 
  405790  & . & 1 & 1 & . & . & . & . & . & 1 & . & . & . & . & . & . & . & . & . \\ 
  811580  & . & . & . & . & . & . & . & . & . & . & . & . & . & . & . & . & 1 & . \\ 
  869550  & . & . & . & . & . & . & . & . & 1 & . & 1 & . & . & . & . & . & . & 1 \\ 
  973896  & . & . & 1 & . & . & 1 & . & . & 1 & . & . & . & . & . & . & . & . & . \\ 
 1014475  & . & 1 & . & . & . & . & . & . & . & . & 1 & 1 & . & . & . & . & . & . \\ 
 1014475  & . & . & . & . & . & . & . & . & . & . & . & . & . & . & . & . & . & 1 \\ 
 1298528  & . & . & . & . & . & . & 1 & . & . & . & . & . & . & . & . & . & . & . \\ 
 1623160  & . & . & . & . & . & . & . & . & . & . & . & . & . & 1 & . & . & 1 & . \\ 
 1623160  & . & . & . & . & . & . & . & . & . & . & . & . & 1 & . & 1 & . & . & . \\ 
 2213400  & . & . & . & . & . & . & . & . & . & 1 & . & . & . & . & 1 & . & . & . \\ 
 3246320  & . & . & . & . & . & . & . & . & . & . & . & . & 1 & . & a & . & . & . \\ 
 3478200  & . & 1 & . & . & . & . & . & . & . & . & 1 & . & 1 & 1 & . & 1 & . & . \\ 
 8115800  & . & . & . & . & . & . & . & . & . & . & . & . & . & . & . & . & 1 & 1 \\ 
18260550  & . & . & . & . & . & . & . & . & . & . & . & . & . & . & . & . & 2 & . 
            \rule[ 0pt]{0pt}{ 4pt} \\ \hline\hline
\end{array}
$$

$a \in \{ 0, 1 \}$
}
\end{table}

\addtocounter{table}{-1}
\begin{table}[h]
\caption{The decomposition matrix modulo $3$ for $B_1$ (cont)}
$$
{\footnotesize
\begin{array}{r|rrrrrrrrrrrrrrrrrr} \hline\hline
\mbox{PIM}& * & * & * & * & * & * & * & * & * & * & * & * & * & * & * & * & * & * \\ \hline
          & 1 & 2 & 2 & 2 & 2 & 2 & 2 & 2 & 2 & 2 & 2 & 3 & 3 & 3 & 3 & 3 & 3 & 3 \\
\Phi      & 9 & 0 & 1 & 2 & 3 & 4 & 5 & 6 & 7 & 8 & 9 & 0 & 1 & 2 & 3 & 4 & 5 & 6 \rule[- 4pt]{0pt}{ 14pt} \\ \hline     
        1 & . & . & . & . & . & . & . & . & . & . & . & . & . & . & . & . & . & . \\ 
      155 & . & . & . & . & . & . & . & . & . & . & . & . & . & . & . & . & . & . \\ 
      187 & . & . & . & . & . & . & . & . & . & . & . & . & . & . & . & . & . & . \\ 
      527 & . & . & . & . & . & . & . & . & . & . & . & . & . & . & . & . & . & . \\ 
    11594 & . & . & . & . & . & . & . & . & . & . & . & . & . & . & . & . & . & . \\ 
    11594 & . & . & . & . & . & . & . & . & . & . & . & . & . & . & . & . & . & . \\ 
    11935 & . & . & . & . & . & . & . & . & . & . & . & . & . & . & . & . & . & . \\ 
    17050 & . & . & . & . & . & . & . & . & . & . & . & . & . & . & . & . & . & . \\ 
    17391 & . & . & . & . & . & . & . & . & . & . & . & . & . & . & . & . & . & . \\ 
    18445 & . & . & . & . & . & . & . & . & . & . & . & . & . & . & . & . & . & . \\ 
    28985 & . & . & . & . & . & . & . & . & . & . & . & . & . & . & . & . & . & . \\ 
    29512 & . & . & . & . & . & . & . & . & . & . & . & . & . & . & . & . & . & . \\ 
    46376 & . & . & . & . & . & . & . & . & . & . & . & . & . & . & . & . & . & . \\ 
    81158 & . & . & . & . & . & . & . & . & . & . & . & . & . & . & . & . & . & . \\ 
   173910 & . & . & . & . & . & . & . & . & . & . & . & . & . & . & . & . & . & . \\ 
   202895 & . & . & . & . & . & . & . & . & . & . & . & . & . & . & . & . & . & . \\ 
   202895 & . & . & . & . & . & . & . & . & . & . & . & . & . & . & . & . & . & . \\ 
   202895 & . & . & . & . & . & . & . & . & . & . & . & . & . & . & . & . & . & . \\ 
   202895 & 1 & . & . & . & . & . & . & . & . & . & . & . & . & . & . & . & . & . \\ 
   371008 & . & 1 & . & . & . & . & . & . & . & . & . & . & . & . & . & . & . & . \\ 
   405790 & . & . & 1 & . & . & . & . & . & . & . & . & . & . & . & . & . & . & . \\ 
   405790 & . & . & . & 1 & . & . & . & . & . & . & . & . & . & . & . & . & . & . \\ 
   405790 & . & . & . & . & 1 & . & . & . & . & . & . & . & . & . & . & . & . & . \\ 
   811580 & . & . & . & . & . & 1 & . & . & . & . & . & . & . & . & . & . & . & . \\ 
   869550 & . & . & . & . & . & . & 1 & . & . & . & . & . & . & . & . & . & . & . \\ 
   973896 & . & . & . & . & . & . & . & 1 & . & . & . & . & . & . & . & . & . & . \\ 
  1014475 & . & . & . & . & . & . & . & . & 1 & . & . & . & . & . & . & . & . & . \\ 
  1014475 & . & . & . & . & . & . & . & . & . & 1 & . & . & . & . & . & . & . & . \\ 
  1298528 & . & . & . & . & . & . & . & . & . & . & 1 & . & . & . & . & . & . & . \\ 
  1623160 & . & . & . & . & . & . & . & . & . & . & . & 1 & . & . & . & . & . & . \\ 
  1623160 & 1 & . & . & . & . & . & . & . & . & . & . & . & 1 & . & . & . & . & . \\ 
  2213400 & . & . & . & . & . & 1 & . & . & . & . & . & . & . & 1 & . & . & . & . \\ 
  3246320 & . & . & . & 1 & . & . & . & . & . & . & . & . & . & . & 1 & . & . & . \\ 
  3478200 & . & . & 1 & . & . & . & . & . & . & . & . & . & . & . & . & 1 & . & . \\ 
  8115800 & 2 & 1 & . & . & . & . & 1 & . & . & 1 & . & . & . & . & . & . & 1 & . \\ 
 18260550 & . & . & . & . & . & . & . & . & . & . & 1 & . & . & . & . & . & 1 & 1 
            \rule[ 0pt]{0pt}{ 4pt} \\ \hline\hline
\end{array}
}
$$
\end{table}

\clearpage

%\addtocounter{table}{-1}
\begin{table}[h]
\caption{\label{Mod3B2PS1} The first basic set of projectives and
some further projectives for $B_2$}
$$
\begin{array}{r|rrrrrrrrrr|rrrr} \hline\hline
\Psi & 1 & 2 & 3 & 4 & 5 & 6 & 7 & 8 & 9 & 10 & 11 & 12 & 13 & 14 \rule[- 4pt]{0pt}{ 14pt} \\ \hline     
      495 & 1 & . & . & . & . & . & . & . & . & . & 1 & . & . & . \\
     6138 & . & 1 & . & . & . & . & . & . & . & . & . & . & . & . \\
    40392 & 1 & . & 1 & . & . & . & . & . & . & . & 1 & 1 & . & . \\
   313038 & . & . & 1 & 1 & . & . & . & . & . & . & 1 & 1 & . & . \\
   417384 & 1 & . & . & . & 1 & . & . & . & . & . & 3 & 1 & . & . \\
   515592 & 1 & . & . & . & . & 1 & . & . & . & . & . & . & 1 & . \\
   782595 & . & 1 & . & . & . & . & 1 & . & . & . & . & . & 3 & . \\
  1826055 & . & . & . & . & . & . & . & 1 & . & . & . & . & . & 3 \\
  3142656 & . & 1 & . & . & . & . & . & . & 1 & . & . & . & . & 1 \\
 11784960 & . & . & 2 & 1 & . & . & 1 & 2 & . & 1 & 1 & 1 & . & 3 
            \rule[ 0pt]{0pt}{ 10pt} \\ \hline\hline
\end{array}
$$
\end{table}

\begin{table}[h]
\caption{\label{Mod3B2PS2} The decomposition matrix modulo $3$ and some relations 
for $B_2$}
$$
\begin{array}{r|rrrrrrrrrr|rrrrrrr} \hline\hline
%\Psi & 1 & 2 & 3 & 4 & 5 & 6 & 7 & 8 & 9 & 10 & 18 & 66 & 106 & 153 & 169 & \otimes \rule[- 4pt]{0pt}{ 10pt} \\ \hline     
\mbox{PIM}& * & * & * & * & * & * & * & * & * & * &   &   &   &   &   &   & \\ \hline
     \Phi & 1 & 2 & 3 & 4 & 5 & 6 & 7 & 8 & 9 & 10& 4 &10 &11 &16 &20 &23 & \otimes_{2,1} \rule[- 4pt]{0pt}{ 10pt} \\ \hline     
      495 & 1 & . & . & . & . & . & . & . & . & . &-1 & 1 & 1 & . & 1 & . &-1 \\
     6138 & . & 1 & . & . & . & . & . & . & . & . & 1 &-1 & . &-1 &-1 &-2 & . \\
    40392 & 1 & . & 1 & . & . & . & . & . & . & . & 1 & . &-1 & . &-1 & . & . \\
   313038 & . & . & 1 & 1 & . & . & . & . & . & . & . & . & 1 & . & . &-1 & . \\
   417384 & 1 & . & . & . & 1 & . & . & . & . & . & . & . & 1 & . & . & . & 1 \\
   515592 & . & . & . & . & . & 1 & . & . & . & . & . & 1 & . & 1 & . & 1 & . \\
   782595 & . & 1 & . & . & . & . & 1 & . & . & . & . & 1 & 1 & . & 1 & . & . \\
  1826055 & . & . & . & . & . & . & . & 1 & . & . & . & . & . & . & . &-1 & . \\
  3142656 & . & 1 & . & . & . & . & . & . & 1 & . & . & . & . & 1 & . & 2 & . \\
 11784960 & . & . & 1 & 1 & . & . & . & 1 & . & 1 & . & . & . & . & 1 & 1 & .  
            \rule[ 0pt]{0pt}{ 10pt} \\ \hline\hline
\end{array}
$$
\end{table}

%relations[2]{[4,10,11,16,20,23]};
%[ [ -1,  1,  1,  0, 0, 0, 0,  0, 0, 0 ],
%  [  1, -1,  0,  0, 0, 1, 1,  0, 0, 0 ],
%  [  1,  0, -1,  1, 1, 0, 1,  0, 0, 0 ]
%  [  0, -1,  0,  0, 0, 1, 0,  0, 1, 0 ],
%  [  1, -1, -1,  0, 0, 0, 1,  0, 0, 1 ],
%  [  0, -2,  0, -1, 0, 1, 0, -1, 2, 1 ] ]

% 155 \otimes 237181 : [ [ -1, 0, 0, 0, 1, 0, 0, 0, 0, 0 ] ]

%blnrs[2]{[4,10,11,16,20,23]};
% [ 18, 66, 71, 106, 153, 169 ]

\begin{table}[h]
\caption{\label{Mod3B3PS1} The first basic set of projectives, a further projective
and some relations for $B_3$}
$$
\begin{array}{r|rrrrrrrrrr|r|rrrr} \hline\hline
\Psi      & 1 & 2 & 3 & 4 & 5 & 6 & 7 & 8 & 9 &10 &11 &15 &22 &20 & \otimes_{3,1}\rule[- 4pt]{0pt}{ 14pt} \\ \hline\hline     
    52173 & 1 & . & . & . & . & . & . & . & . & . & . & . & 1 & . &-1 \\
   365211 & . & 1 & . & . & . & . & . & . & . & . & . &-1 & . &-1 & . \\
   782595 & 1 & . & 1 & . & . & . & . & . & . & . & . & . &-1 & . & . \\
  1095633 & . & . & . & 1 & . & . & . & . & . & . & . & . & . & 1 & . \\
  1826055 & . & . & . & . & 1 & . & . & . & . & . & 3 & 1 & . &-1 & . \\
  1826055 & 1 & . & . & . & . & 1 & . & . & . & . & . & 1 & . & . & 1 \\
  2921688 & . & 1 & . & . & . & . & 1 & . & . & . & 1 & . &-1 & 1 & . \\
  5478165 & . & 1 & . & . & . & . & . & 1 & . & . & . & 1 & 1 & . & . \\
  5478165 & . & . & 1 & . & 1 & . & . & . & 1 & . & . & . & 1 & . & . \\
 14608440 & . & . & . & . & 3 & . & 1 & . & . & 1 & 4 & . & 1 & 1 & .  
            \rule[ 0pt]{0pt}{ 10pt} \\ \hline\hline
\end{array}
$$
\end{table}

% relation 15: [ 0, -1, 0, 0, 1, 1, 0, 1, 0, 0 ]
% relation 20: [ 0, -1, 0, 1, -1, 0, 1, 0, 0, 1 ]
% relation 22: [ 1, 0, -1, 0, 0, 0, -1, 1, 1, 1 ] 
% 155 \otimes 237181: 16337 : [ [ -1, 0, 0, 0, 0, 1, 0, 0, 0, 0 ] ]
% blnrs[3]{[15,20,22]};
% [ 136, 168, 178 ]

\begin{table}[h]
\caption{\label{Mod3B3PS2} The decomposition matrix modulo $3$ for $B_3$}
$$
\begin{array}{r|rrrrrrrrrr} \hline\hline
\mbox{PIM}& * & * & * & * & * & * & * & * & * & * \\ \hline
\Phi      & 1 & 2 & 3 & 4 & 5 & 6 & 7 & 8 & 9 & 10  \rule[- 4pt]{0pt}{ 14pt} \\ \hline\hline     
    52173 & 1 & . & . & . & . & . & . & . & . & .  \\
   365211 & . & 1 & . & . & . & . & . & . & . & .  \\
   782595 & 1 & . & 1 & . & . & . & . & . & . & .  \\
  1095633 & . & . & . & 1 & . & . & . & . & . & .  \\
  1826055 & . & . & . & . & 1 & . & . & . & . & .  \\
  1826055 & 1 & . & . & . & . & 1 & . & . & . & .  \\
  2921688 & . & 1 & . & . & . & . & 1 & . & . & .  \\
  5478165 & . & 1 & . & . & . & . & . & 1 & . & .  \\
  5478165 & . & . & 1 & . & . & . & . & . & 1 & .  \\
 14608440 & . & . & . & . & 1 & . & 1 & . & . & 1   
            \rule[ 0pt]{0pt}{ 10pt} \\ \hline\hline

%
% relation 168: [ 0, -1, 0, 1, -1, 0, 1, 0, 0, 1 ]
% relation 178: [ 1, 0, -1, 0, 0, 0, -1, 1, 1, 1 ] 
\end{array}
$$
\end{table}

\begin{table}[h]
\caption{\label{Mod3B5PS1} The decomposition matrix modulo $3$ 
and a relation for $B_5$}
$$
\begin{array}{r|rrrrr|r} \hline\hline
\mbox{PIM}& * & * & * & * & * &  \\ \hline
     \Phi & 1 & 2 & 3 & 4 & 5 & 8 \rule[- 4pt]{0pt}{ 14pt} \\ \hline\hline
  2347785 & 1 & . & . & . & . & 1 \\
  2347785 & . & 1 & . & . & . & -1 \\
 11738925 & 1 & . & 1 & . & . & -1 \\
 11738925 & . & 1 & . & 1 & . & 1 \\
 23477850 & . & . & 1 & . & 1 & 1  
            \rule[ 0pt]{0pt}{ 10pt} \\ \hline\hline
\end{array}
$$
\end{table}

\begin{table}[h]
\caption{\label{Mod3B6PS1} The decomposition matrix modulo $3$ for $B_6$}
$$
\begin{array}{r|rrrr} \hline\hline
\mbox{PIM}& * & * & * & * \\ \hline
\Phi & 1 & 2 & 3 & 4 \rule[- 4pt]{0pt}{ 14pt} \\ \hline\hline
  3286899 & 1 & . & . & . \\
  3286899 & . & 1 & . & . \\
 13147596 & . & . & 1 & . \\
 13147596 & . & . & . & 1 
            \rule[ 0pt]{0pt}{ 10pt} \\ \hline\hline
\end{array}
$$
\end{table}

%\subsection{Origin of projectives}

\begin{table}[h]
\caption{\label{Mod3B1PS1src} The origin of the first basic set of projectives 
for $B_1$}
$
\begin{array}{c|c} \hline\hline
\Psi & \text{\rm src} \rule[- 4pt]{0pt}{ 14pt} \\ \hline
 1 &  1^- \rule[  0pt]{0pt}{ 12pt} \\
 2 &  4^-  \\
 3 &  6^+  \\
 4 &  3^+  \\
 5 & 21^-  \\
 6 & 18^-  \\
 7 & 17^+  \\
 8 & 21^+  \\
 9 & 23^- 
         \rule[ 0pt]{0pt}{ 10pt} \\ \hline\hline
\end{array}
$\,\,\,
$
\begin{array}{c|c} \hline\hline
\Psi & \text{\rm src} \rule[- 4pt]{0pt}{ 14pt} \\ \hline
10 & 26^+ \rule[  0pt]{0pt}{ 12pt} \\
11 & ( 15^+ + 51^- )/2  \\
12 & 22^-  \\
13 & 24^+  \\
14 & 29^+  \\
15 & 32^+  \\
16 & 31^+ / 2  \\
17 & 37^+  \\
18 & 36^+  
         \rule[ 0pt]{0pt}{ 10pt} \\ \hline\hline
\end{array}
$\,\,\,
$
\begin{array}{c|c} \hline\hline
\Psi & \text{\rm src} \rule[- 4pt]{0pt}{ 14pt} \\ \hline
19 & 36^- \rule[  0pt]{0pt}{ 12pt} \\
20 & 44^+  \\
21 & 53^+  \\
22 & 48^-  \\
23 & 22^+ / 3 \\
24 & 46^+  \\
25 & 52^+  \\
26 & 51^+  \\
27 & 43^-  
         \rule[ 0pt]{0pt}{ 10pt} \\ \hline\hline
\end{array}
$\,\,\,
$
\begin{array}{c|c} \hline\hline
\Psi & \text{\rm src} \rule[- 4pt]{0pt}{ 14pt} \\ \hline
28 & 55^- \rule[  0pt]{0pt}{ 12pt} \\
29 & 59^-  \\
30 & 56^+  \\
31 & 54^+  \\
32 & 56^-  \\
33 & 57^+  \\
34 & 51^-  \\
35 & 49^+  \\
36 & 50^+  \\
11' & 15^+
         \rule[ 0pt]{0pt}{ 10pt} \\ \hline\hline
\end{array}
$
\end{table}

% curPS1[1] := [ indm2[1], indm2[4], indm1[6], indm1[3], indm2[21], indm2[18],
%                indm1[17], indm1[21], indm2[23], indm1[26], indm1[15], indm2[22],
%                indm1[24], indm1[29], indm1[32], indm1[31]/2, indm1[37], indm1[36],
%                indm2[36], indm1[44], indm1[53], indm2[48], indm1[22]/3, indm1[46],
%                indm1[52], indm1[51], indm2[43], indm2[55], indm2[59], indm1[56],
%                indm1[54], indm2[56], indm1[57], indm2[51], indm1[49], indm1[50] ];;

\begin{table}[h]
\caption{\label{Mod3B1Furthsrc} The origin of the further projectives for 
$B_1$}
$
\begin{array}{c|c} \hline\hline
\Psi & \text{\rm src} \rule[- 4pt]{0pt}{ 14pt} \\ \hline
37 & 47^+ \rule[  0pt]{0pt}{ 12pt} \\
38 & 45^-  \\
39 & 40^-  \\
40 & 46^-  \\
41 & 34^+  \\
42 & 27^-/2 
         \rule[ 0pt]{0pt}{ 10pt} \\ \hline\hline
\end{array}
$\,\,\,
% [ [ 47, "+" ], [ 45, "-" ], [ 40, "-" ], [ 46, "-" ], [ 34, "+" ], [ 27, "-" ] ]
$
\begin{array}{c|c} \hline\hline
\Psi & \text{\rm src} \rule[- 4pt]{0pt}{ 14pt} \\ \hline
43 & 25^+ \rule[  0pt]{0pt}{ 12pt} \\
44 & 30^+ \\
45 & 28^-  \\
46 & 19^+  \\
47 & 20^+  \\
48 & 11^+  
         \rule[ 0pt]{0pt}{ 10pt} \\ \hline\hline
\end{array}
$\,\,\,
% [ [ 25, "+" ], [ 30, "+" ], [ 28, "-" ], [ 19, "+" ], [ 20, "+" ], [ 11, "+" ] ]
$
\begin{array}{c|c} \hline\hline
\Psi & \text{\rm src} \rule[- 4pt]{0pt}{ 14pt} \\ \hline
49 &  8^- \rule[  0pt]{0pt}{ 12pt} \\
50 &  9^+  \\
51 &  7^+  \\
52 &  7^-  \\
53 & 14^+  \\
54 &  3^- 
         \rule[ 0pt]{0pt}{ 10pt} \\ \hline\hline
\end{array}
$\,\,\,
% [ [ 8, "-" ], [ 9, "+" ], [ 7, "+" ], [ 7, "-" ], [ 14, "+" ], [ 3, "-" ] ]
$
\begin{array}{c|c} \hline\hline
\Psi & \text{\rm src} \rule[- 4pt]{0pt}{ 14pt} \\ \hline
55 &  5^- \rule[  0pt]{0pt}{ 12pt} \\
56 &  5^+  \\
57 &  6^-  \\
68 &  2^+  \\
59 &  2^-  \\
60 &  4^+  \\
61 &  1^+\
         \rule[ 0pt]{0pt}{ 10pt} \\ \hline\hline
\end{array}
$
% [ [ 5, "-" ], [ 5, "+" ], [ 6, "-" ], [ 2, "+" ], [ 2, "-" ], [ 4, "+" ], [ 1, "+" ] ]

\end{table}

% proofProj5OriginOriginal :=
% [ 47, 111, 106, 112, 34, 93, 25, 30, 94, 19, 20, 11, 74, 9, 7, 73, 14, 69, 71, 5, 72, 2, 68, 4, 1 ];
% [1..66]: indM1 (+), [67..130]: indM2 (-)

\begin{table}[h]
\caption{\label{Mod3B2356PS1src} The origin of the relevant projectives for $B_2$,~$B_3$,~$B_5$ and~$B_6$}
$
\begin{array}{c|c} \hline\hline
\multicolumn{2}{c}{B_2} \rule[- 4pt]{0pt}{ 16pt} \\ \hline\hline
\Psi & \text{\rm src} \rule[- 4pt]{0pt}{ 16pt} \\ \hline     
 1 &  1^- / 3 \rule[  0pt]{0pt}{ 12pt} \\
 2 &  4^- / 3  \\
 3 &  6^+ / 3 \\
 4 &  9^+ / 3 \\
 5 & 21^+ / 4 \\
 6 & 17^+ / 3 \\
 7 & 11^+ / 3 \\
 8 & 14^- / 4 \\
 9 & 18^+ \\
10 & 50^+ \\
12 &  5^+ \\
13 & 10^+ \\
14 & 39^+ \\
15 & 42^+ 
         \rule[ 0pt]{0pt}{ 10pt} \\ \hline\hline
\end{array}
$\,\,\,
%\end{table}
% curPS1[2] := [ indm2[1]/3, indm2[4]/3, indm1[6]/3, indm1[9]/3, indm1[21]/4,
%                indm1[17]/3, indm1[11]/3, indm2[14]/4, indm1[18], indm1[50],
%                indm1[5], indm1[10], indm1[39], indm1[42] ];;
$
\begin{array}{c|c} \hline\hline
\multicolumn{2}{c}{B_3} \rule[- 4pt]{0pt}{ 16pt} \\ \hline\hline
\Psi & \text{\rm src} \rule[- 4pt]{0pt}{ 16pt} \\ \hline     
 1 & 10^+ \rule[  0pt]{0pt}{ 12pt} \\
 2 & 47^- \\
 3 &  5^+ \\
 4 & 57^+ \\
 5 & 55^+ \\
 6 & 22^+ / 2 \\
 7 & 56^+ \\
 8 & 32^+ / 7 \\
 9 & 50^+ / 2 \\
10 & 49^+ / 2 \\
11 & 41^+ 
         \rule[ 0pt]{0pt}{ 10pt} \\ \hline\hline
\end{array}
$\,\,\,
% curPS1[3] := [ indm1[10], indm2[47], indm1[5], indm1[57], indm1[55],
%                indm1[22]/2, indm1[56], indm1[32]/7, indm1[50]/2, indm1[49]/2, indm1[41] ];;
$
\begin{array}{c|c} \hline\hline
\multicolumn{2}{c}{B_5} \rule[- 4pt]{0pt}{ 16pt} \\ \hline\hline
\Psi & \text{\rm src} \rule[- 4pt]{0pt}{ 16pt} \\ \hline     
 1 & 3^+ / 2 \rule[  0pt]{0pt}{ 12pt} \\
 2 & 57^+ \\
 3 & 31^+ / 7 \\
 4 &  2^+ \\
 5 & 49^+
         \rule[ 0pt]{0pt}{ 10pt} \\ \hline\hline
\end{array}
$\,\,\,
$
\begin{array}{c|c} \hline\hline
\multicolumn{2}{c}{B_6} \rule[- 4pt]{0pt}{ 16pt} \\ \hline\hline
\Psi & \text{\rm src} \rule[- 4pt]{0pt}{ 16pt} \\ \hline     
 1 & 17^+ \rule[  0pt]{0pt}{ 12pt} \\
 2 & 44^+ \\
 3 &  1^+ \\
 4 &  2^+
         \rule[ 0pt]{0pt}{ 10pt} \\ \hline\hline
\end{array}
$
\end{table}

%\begin{table}[h]
%\caption{\label{Mod3B3PS1src} The origin of the first basic set of projectives for Block~$3$}
%$$
%\begin{array}{c|c} \hline\hline
%\Psi & \text{\rm src} \rule[- 4pt]{0pt}{ 14pt} \\ \hline     
% 1 & 10^+ \rule[  0pt]{0pt}{ 12pt} \\
% 2 & 47^- \\
% 3 &  5^+ \\
% 4 & 57^+ \\
% 5 & 55^+ \\
% 6 & 22^+ / 2 \\
% 7 & 56^+ \\
% 8 & 32^+ / 7 \\
% 9 & 50^+ / 2 \\
%10 & 49^+ / 2 \\
%11 & 41^+ 
%         \rule[ 0pt]{0pt}{ 10pt} \\ \hline\hline
%\end{array}
%% curPS1[3] := [ indm1[10], indm2[47], indm1[5], indm1[57], indm1[55],
%%                indm1[22]/2, indm1[56], indm1[32]/7, indm1[50]/2, indm1[49]/2, indm1[41] ];;
%$$
%\end{table}
%
%\begin{table}[h]
%\caption{\label{Mod3B5PS1src} The origin of the first basic set of projectives for Block~$5$}
%$$
%\begin{array}{c|c} \hline\hline
%\Psi & \text{\rm src} \rule[- 4pt]{0pt}{ 14pt} \\ \hline     
% 1 & 3^+ / 2 \rule[  0pt]{0pt}{ 12pt} \\
% 2 & 57^+ \\
% 3 & 31^+ / 7 \\
% 4 &  2^+ \\
% 5 & 49^+
%         \rule[ 0pt]{0pt}{ 10pt} \\ \hline\hline
%\end{array}
%$$
%\end{table}
%
%\begin{table}[h]
%\caption{\label{Mod3B6PS1src} The origin of the first basic set of projectives for Block~$6$}
%$$
%\begin{array}{c|c} \hline\hline
%\Psi & \text{\rm src} \rule[- 4pt]{0pt}{ 14pt} \\ \hline     
% 1 & 17^+ \rule[  0pt]{0pt}{ 12pt} \\
% 2 & 44^+ \\
% 3 &  1^+ \\
% 4 &  2^+
%         \rule[ 0pt]{0pt}{ 10pt} \\ \hline\hline
%\end{array}
%$$
%\end{table}

\begin{table}[h]
\caption{\label{Mod3B1Rels} Some relations in $B_1$} 
$$
{\footnotesize
\hspace*{-15.00mm}
\begin{array}{r|rrrrrrrrrrrrrrrrrr} \hline\hline
 \text{rel} & 4 & 14 & 17 & 18 & 19 & 23 & 29 & 30 & 33 & 40 & 49 & 51 & 52 & 54 & 56 & 59 & 60 & 65 \\ \hline\hline
 1 & -1 & . & 1 & 1 & 1 & -1 & 1 & 1 & -1 & . & -1 & . & 1 & -1 & -1 & . & -1 & . \rule[- 4pt]{0pt}{ 14pt} \\ 
 155 & 1 & . & -1 & -1 & -1 & 1 & . & . & . & . & 1 & . & . & 1 & 1 & . & 1 & . \\ 
 187 & 1 & . & . & -1 & -1 & 1 & -1 & -1 & 1 & . & 1 & . & -1 & 1 & 1 & . & 1 & . \\ 
 527 & . & -1 & . & . & . & 1 & -1 & -1 & 1 & . & 1 & . & . & . & . & . & . & 1 \\ 
 11594 & . & . & 1 & 1 & . & . & 1 & . & . & 1 & . & . & 1 & . & . & . & . & . \\ 
 11594 & . & 1 & . & . & . & -1 & . & 1 & -1 & . & -1 & . & . & . & . & . & . & -1 \\ 
 11935 & . & . & 1 & . & 1 & -1 & . & . & -1 & . & -1 & . & . & . & -1 & . & -1 & . \\ 
 17050 & . & . & . & 1 & . & . & 1 & 1 & . & . & . & . & 1 & . & . & . & . & . \\ 
 17391 & . & . & . & . & 1 & . & . & . & . & -1 & . & . & . & -1 & -1 & . & -1 & . \\ 
 18445 & . & . & . & . & . & -1 & 1 & 1 & -1 & . & . & . & 1 & . & . & . & . & . \\ 
 28985 & . & . & 1 & 1 & 1 & . & 1 & 1 & . & . & . & . & 1 & -1 & -1 & . & -1 & . \\ 
 29512 & . & 1 & . & . & . & -1 & . & . & . & . & -1 & . & -1 & . & . & . & . & -1 \\ 
 46376 & . & . & . & . & . & . & -1 & -1 & . & -1 & . & . & -1 & . & -1 & . & -2 & . \\ 
 81158 & . & . & . & . & . & 1 & . & . & 1 & . & 1 & . & 1 & . & 1 & . & 2 & 1 \\ 
 173910 & . & . & . & . & . & 1 & . & . & 1 & . & 1 & . & . & . & . & . & . & . \\ 
 202895 & . & . & . & . & . & . & 1 & . & . & . & . & 1 & 1 & . & . & . & . & . \\ 
 202895 & . & . & . & . & . & . & . & . & -1 & . & -1 & 1 & . & . & -1 & -1 & -2 & . \\ 
 202895 & . & . & . & . & . & . & . & 1 & . & . & 1 & -1 & . & . & . & . & -1 & . \\ 
 202895 & . & . & . & . & . & . & . & . & . & 1 & . & . & . & 1 & 1 & . & 2 & . \\ 
 371008 & . & . & . & . & . & . & . & . & 1 & . & . & . & . & . & 1 & . & 2 & . \\ 
 405790 & . & . & . & . & . & . & . & . & . & . & . & . & . & . & . & 1 & . & . \\ 
 405790 & . & . & . & . & . & . & . & . & . & . & . & . & . & . & 1 & . & 1 & . \\ 
 405790 & . & . & . & . & . & . & . & . & . & 1 & . & . & . & 1 & . & . & . & . \\ 
 811580 & . & . & . & . & . & . & . & . & . & . & 1 & . & . & . & . & . & . & . \\ 
 869550 & . & . & . & . & . & . & . & . & . & . & . & . & . & 1 & 1 & . & 1 & . \\ 
 973896 & . & . & . & . & . & . & . & . & . & . & . & . & . & . & . & . & . & 1 \\ 
 1014475 & . & . & . & . & . & . & . & . & . & . & . & . & 1 & . & . & . & . & 1 \\ 
 1014475 & . & . & . & . & . & . & . & . & . & . & . & 1 & . & . & . & . & . & . \\ 
 1298528 & . & . & . & . & . & . & . & . & . & . & . & . & . & . & . & . & . & . \\ 
 1623160 & . & . & . & . & . & . & . & . & . & . & . & . & . & . & . & 1 & . & . \\ 
 1623160 & . & . & . & . & . & . & . & . & . & . & . & . & . & . & . & . & . & . \\ 
 2213400 & . & . & . & . & . & . & . & . & . & . & . & . & . & . & . & . & . & . \\ 
 3246320 & . & . & . & . & . & . & . & . & . & . & . & . & . & . & . & . & . & . \\ 
 3478200 & . & . & . & . & . & . & . & . & . & . & . & . & . & . & . & . & . & . \\ 
 8115800 & . & . & . & . & . & . & . & . & . & . & . & . & . & . & . & . & . & . \\ 
 18260550 & . & . & . & . & . & . & . & . & . & . & . & . & . & . & . & . & . & . 
            \rule[ 0pt]{0pt}{ 4pt} \\ \hline\hline
\end{array}
}
$$
\end{table}

\clearpage

\begin{table}[h]
\caption{Some relations in $B_1$ (cont)} 
$$
{\footnotesize
\hspace*{-15.00mm}
\begin{array}{r|rrrrrrrrrrrrrrrrrr} \hline\hline
\text{rel} &  66 & 70 & 71 & 74 & 77 & 79 & 80 & 82 & 84 & 88 & 98 & 104 & 105 & 110 & 112 & \otimes_{1,1} & \otimes_{1,2} & \otimes_{1,3}  \rule[- 4pt]{0pt}{ 14pt} \\ \hline\hline
 1 & . & 1 & 1 & 1 & . & 1 & . & . & 1 & 1 & 1 & 1 & 1 & 1 & 1 & 1 & . & -2 \\ 
 155 & . & -1 & -1 & -1 & . & -1 & . & . & -1 & -1 & -1 & -1 & -1 & -1 & -1 & -1 & 1 & 2 \\ 
 187 & . & -1 & -1 & -1 & . & -1 & . & . & -1 & . & -1 & -1 & -1 & -1 & -1 & -2 & . & 2 \\ 
 527 & . & -1 & -1 & -1 & . & -1 & . & . & -1 & . & . & -1 & . & -2 & -2 & . & . & . \\ 
 11594 & . & . & . & . & . & . & . & . & . & . & . & . & . & . & . & . & . & . \\ 
 11594 & . & 1 & 1 & 1 & . & 1 & . & . & 1 & . & . & 1 & . & 2 & 2 & -1 & . & . \\ 
 11935 & . & 1 & 1 & 1 & . & 1 & . & . & 1 & 1 & 1 & 1 & 1 & 1 & 2 & . & . & -3 \\ 
 17050 & . & . & . & . & . & . & . & . & . & . & . & . & . & . & . & 2 & . & . \\ 
 17391 & . & . & -1 & . & . & -1 & . & -1 & . & . & 1 & . & 1 & -1 & . & 1 & . & -1 \\ 
 18445 & . & 1 & 1 & 1 & . & 1 & . & 1 & 1 & . & . & 1 & . & 1 & 1 & 1 & . & . \\ 
 28985 & . & . & . & . & -1 & -1 & -1 & -1 & . & . & 1 & . & . & -1 & -1 & 1 & -1 & -1 \\ 
 29512 & . & 1 & 1 & 1 & . & 1 & . & . & 1 & . & . & 1 & . & 2 & 2 & . & -1 & . \\ 
 46376 & -1 & . & . & . & . & . & . & . & . & . & 1 & 1 & 1 & -1 & . & . & . & -3 \\ 
 81158 & . & -1 & -1 & -1 & . & -1 & -1 & . & -1 & -1 & -2 & -1 & -2 & . & -1 & . & . & 4 \\ 
 173910 & . & -1 & -1 & -1 & . & -1 & . & -1 & -1 & . & . & -1 & . & -2 & -2 & . & . & . \\ 
 202895 & . & 1 & . & . & . & . & -1 & . & . & -1 & . & . & -1 & . & . & . & . & . \\ 
 202895 & . & 1 & 1 & 1 & . & 1 & . & . & 1 & . & 1 & 1 & 1 & 1 & 1 & . & . & -4 \\ 
 202895 & . & -1 & . & . & -1 & -1 & . & . & . & . & 1 & . & . & -1 & -1 & . & . & 1 \\ 
 202895 & . & . & . & -1 & 1 & 1 & . & . & -1 & . & -2 & -1 & -1 & . & . & . & . & 2 \\ 
 371008 & 1 & . & . & . & . & . & . & . & . & . & -1 & -1 & -1 & 1 & . & . & . & 2 \\ 
 405790 & . & . & 1 & 1 & -1 & . & . & . & . & . & 1 & 1 & . & . & . & . & . & . \\ 
 405790 & 1 & . & . & . & . & . & 1 & . & -1 & . & . & . & . & -1 & . & . & . & . \\ 
 405790 & . & . & . & . & . & 1 & . & . & . & . & . & . & . & . & . & . & . & . \\ 
 811580 & . & -1 & -1 & -1 & . & -1 & . & . & . & . & . & -1 & . & -1 & . & . & . & 1 \\ 
 869550 & . & . & . & . & . & 1 & . & 1 & . & . & -1 & . & -1 & 1 & 1 & . & . & 1 \\ 
 973896 & . & . & 1 & . & . & 1 & . & 1 & . & . & . & 1 & . & . & . & 1 & . & . \\ 
 1014475 & . & . & . & . & . & . & . & 1 & . & . & . & . & . & . & . & . & 1 & . \\ 
 1014475 & . & 1 & . & . & 1 & 1 & . & . & . & . & -1 & . & . & 1 & 1 & . & . & 1 \\ 
 1298528 & 1 & . & . & . & . & . & 1 & . & . & 1 & 1 & . & 1 & -1 & -1 & . & . & -1 \\ 
 1623160 & . & . & . & . & . & . & . & . & 1 & 1 & 1 & . & 1 & . & . & . & . & . \\ 
 1623160 & . & . & . & 1 & . & . & . & . & 1 & . & . & . & . & 2 & 1 & . & . & . \\ 
 2213400 & . & 1 & 1 & 1 & . & 1 & . & 1 & . & . & . & 1 & . & 1 & 1 & . & . & . \\ 
 3246320 & . & . & . & . & . & . & . & . & 1 & . & . & . & . & 1 & . & . & . & . \\ 
 3478200 & . & . & . & . & 1 & . & 1 & . & . & 1 & . & . & 1 & . & . & . & . & . \\ 
 8115800 & . & . & . & . & . & . & . & . & . & . & 1 & 1 & 1 & -1 & . & . & . & -1 \\ 
 18260550 & . & . & . & . & . & . & . & . & . & . & . & . & . & 1 & 1 & . & . & 1 
            \rule[ 0pt]{0pt}{ 4pt} \\ \hline\hline
\end{array}
}
$$
\end{table}

\clearpage

\begin{table}[h]
\caption{\label{Mod5B1PS1} The decomposition matrix modulo $5$ for $B_1$}
$$
\begin{array}{r|rrrrrrrrrrrrrr} \hline\hline
\mbox{PIM}& * & * & * & * & * & * & * & * & * & * & * & * & * & * \\ \hline
          &   &   &   &   &   &   &   &   &   & 1 & 1 & 1 & 1 & 1 \\
\Phi      & 1 & 2 & 3 & 4 & 5 & 6 & 7 & 8 & 9 & 0 & 1 & 2 & 3 & 4 \rule[- 4pt]{0pt}{ 14pt} \\ \hline
        1 & 1 & . & . & . & . & . & . & . & . & . & . & . & . & . \\
     6138 & . & 1 & . & . & . & . & . & . & . & . & . & . & . & . \\
    11594 & 1 & . & 1 & . & . & . & . & . & . & . & . & . & . & . \\
    11594 & 1 & . & . & 1 & . & . & . & . & . & . & . & . & . & . \\
    29512 & . & 1 & . & . & 1 & . & . & . & . & . & . & . & . & . \\
    40392 & . & . & . & 1 & . & 1 & . & . & . & . & . & . & . & . \\
   278256 & . & . & 1 & . & . & . & 1 & . & . & . & . & . & . & . \\
   371008 & . & . & . & . & 1 & . & . & 1 & . & . & . & . & . & . \\
   649264 & 1 & . & . & 1 & . & 1 & . & . & 1 & . & . & . & . & . \\
  2921688 & . & . & . & . & . & . & . & . & . & 1 & . & . & . & . \\
  3286899 & . & . & . & . & . & . & . & 1 & . & . & 1 & . & . & . \\
  3339072 & 1 & . & 1 & . & . & . & . & . & 1 & . & . & 1 & . & . \\
  4382532 & . & 1 & . & . & 1 & 1 & . & 1 & . & . & . & . & 1 & . \\
 11269368 & . & 1 & . & . & . & . & . & . & . & . & . & . & 1 & 1 
            \rule[ 0pt]{0pt}{ 4pt} \\ \hline\hline
\end{array}
$$
\end{table}

\begin{table}[h]
\caption{\label{Mod5B3PS1} The decomposition matrix modulo $5$ for $B_3$}
$$
\begin{array}{r|rrrrrrrrrrrrrr} \hline\hline
\mbox{PIM}& * & * & * & * & * & * & * & * & * & * & * & * & * & * \\ \hline
          &   &   &   &   &   &   &   &   &   & 1 & 1 & 1 & 1 & 1 \\
\Phi      & 1 & 2 & 3 & 4 & 5 & 6 & 7 & 8 & 9 & 0 & 1 & 2 & 3 & 4 \rule[- 4pt]{0pt}{ 14pt} \\ \hline 
      187 & 1 & . & . & . & . & . & . & . & . & . & . & . & . & . \\
      527 & . & 1 & . & . & . & . & . & . & . & . & . & . & . & . \\
    46376 & 1 & 1 & 1 & . & . & . & . & . & . & . & . & . & . & . \\
    46376 & . & . & . & 1 & . & . & . & . & . & . & . & . & . & . \\
    57288 & . & 1 & . & . & 1 & . & . & . & . & . & . & . & . & . \\
   365211 & . & . & . & . & . & 1 & . & . & . & . & . & . & . & . \\
   417384 & . & 1 & 1 & . & 1 & . & 1 & . & . & . & . & . & . & . \\
   556512 & 1 & . & 1 & . & . & . & . & 1 & . & . & . & . & . & . \\
   944384 & . & . & . & 1 & . & . & . & . & 1 & . & . & . & . & . \\
  1292544 & . & . & . & . & . & 1 & . & . & . & 1 & . & . & . & . \\
  3142656 & . & . & . & . & . & . & . & . & 1 & . & 1 & . & . & . \\
  5936128 & . & . & . & . & . & 1 & 1 & . & . & 1 & . & 1 & . & . \\
  5936128 & . & . & 1 & . & . & . & 1 & 1 & . & . & . & . & 1 & . \\
 26295192 & . & . & . & . & . & 1 & . & . & . & . & . & 1 & . & 1 
            \rule[ 0pt]{0pt}{ 4pt} \\ \hline\hline
\end{array}
$$
\end{table}

\begin{table}[h]
\caption{\label{Mod5B4PS1} The decomposition matrix modulo $5$ for $B_4$}
$$
\begin{array}{r|rrrrrrrrrrrrrr} \hline\hline
\mbox{PIM}& * & * & * & * & * & * & * & * & * & * & * & * & * & * \\ \hline
          &   &   &   &   &   &   &   &   &   & 1 & 1 & 1 & 1 & 1 \\
\Phi      & 1 & 2 & 3 & 4 & 5 & 6 & 7 & 8 & 9 & 0 & 1 & 2 & 3 & 4 \rule[- 4pt]{0pt}{ 14pt} \\ \hline 
      341 & 1 & . & . & . & . & . & . & . & . & . & . & . & . & . \\
    17391 & . & 1 & . & . & . & . & . & . & . & . & . & . & . & . \\
    40579 & 1 & . & 1 & . & . & . & . & . & . & . & . & . & . & . \\
    81158 & . & . & . & 1 & . & . & . & . & . & . & . & . & . & . \\
   313038 & 1 & 1 & . & . & 1 & . & . & . & . & . & . & . & . & . \\
   515592 & . & . & . & 1 & . & 1 & . & . & . & . & . & . & . & . \\
   973896 & . & 1 & . & . & 1 & . & 1 & . & . & . & . & . & . & . \\
   973896 & . & . & 1 & . & . & . & . & 1 & . & . & . & . & . & . \\
  1095633 & . & . & . & . & . & . & . & . & 1 & . & . & . & . & . \\
  1298528 & . & . & . & . & . & 1 & . & . & . & 1 & . & . & . & . \\
  1947792 & 1 & . & 1 & . & 1 & . & . & . & . & . & 1 & . & . & . \\
  3286899 & . & . & . & . & . & . & . & . & . & 1 & . & 1 & . & . \\
  9860697 & . & . & . & 1 & . & 1 & 1 & . & . & 1 & . & . & 1 & . \\
 13147596 & . & . & . & . & . & . & . & . & 1 & . & . & 1 & . & 1
            \rule[ 0pt]{0pt}{ 4pt} \\ \hline\hline
\end{array}
$$
\end{table}

\end{document}